\documentstyle[11pt]{article}

\newtheorem{stat}{Statement}
\newcommand{\proof}{{\bf Proof: }}
\newcommand{\proofend}{{\bf Q.E.D.}\vspace{0.5cm}}

\title{
Current Mathematics Appears to Be Inconsistent
}

\author{
\vspace{0.5cm}
Guang-Liang Li 
\footnote{The corresponding author.}
\hspace{2cm}
Victor O. K. Li
\\
Department of Electrical and Electronic Engineering \\
The University of Hong Kong \\
Room 601, Chow Yei Ching Bldg. \\
Pokfulam Road, Hong Kong \\
\{glli,vli\}@eee.hku.hk\\
Phone: (852)2857 8495,
Fax: (852)2559 8738
}

\begin{document}

\maketitle


\begin{abstract}
{\em 
We show that
some mathematical results and their negations
are both deducible. 
The derived contradictions
indicate the inconsistency of current mathematics.
}

This paper is an updated version of 
arXiv:math/0606635v3 with additional results and proofs.

\end{abstract}

According to G{\" o}del's second incompleteness theorem,
if mathematics is consistent, the consistency cannot
be proved by mathematical means \cite{Ebbi}.
However,
to show the inconsistency of mathematics,
it is sufficient to deduce both
a mathematical statement and its negation.
This paper gives some examples of such
statements.
We first consider an example in a geometric setting.

Let line $l$ be the $x$-axis 
on the Euclidean plane.
Let $X_n, n = 1, 2, \cdots$ be points on $l$ with
coordinates $(x_n, 0)$, i.e.,
$X_n = (x_n, 0)$, where
$x_n, n = 1, 2, \cdots$ form
a strictly decreasing sequence of
positive numbers, converging to a unique limit $a$.
Denote by $A$ the point $(a, 0)$ on $l$, and
$B$ the point $(a, b)$ on the plane,
where $b$
is an arbitrarily given positive number (Figure 1).

\begin{stat}
\label{st-G1}
There is a line segment $AA'$
on $l$ from point $A = (a, 0)$ to
point $A' = (a', 0)$, 
where $a' > a$, such that there is not
any point $X_n = (x_n, 0)$ on $AA'$.
\end{stat}

\proof
For each point $X_n = (x_n, 0)$ with $n > 1$, 
there is a line segment $X_nB$ from $X_n$ to
$B$. 
Since $x_n > a$ for
any $x_n$, no $X_nB$ is 
perpendicular to $l$.
Let $X_nBX_1$ represent the set of the points that form
the triangle with sides
$X_nB, BX_1$, and $X_1X_n$. 
Denote by $\Delta$ the union of all such triangles
$X_nBX_1$. 
Let $ABX_1$ represent the set of the points that form
the triangle with sides
$AB, BX_1$, and $X_1A$. 
The line passing through both $A$ and $B$ is
perpendicular to $l$.
No $X_nB$ coincides with $AB$. So
$\Delta$ is a proper subset of
$ABX_1$ (Figure 1). 

\setlength{\unitlength}{0.50mm}
\begin{picture}(250,130)
\put(80,95){\line(2,-1){148}}
\put(80,95){\line(1,-1){74}}
\put(60,20){\vector(1,0){190}}
\put(79,15){\makebox(0,0){$A$}}
\put(55,18){\makebox(0,0){$l$}}
\put(255,18){\makebox(0,0){$x$}}
\put(156,15){\makebox(0,0){$X_n$}}
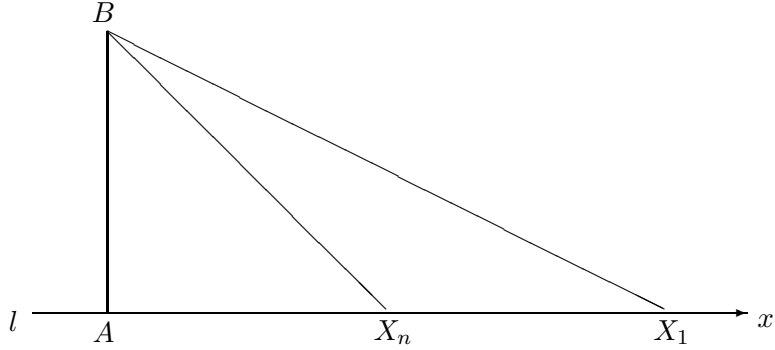
\put(160,0){\makebox(0,0){Figure 1. The geometric setting.}}
\put(230,15){\makebox(0,0){$X_1$}}
\put(80,20){\line(0,1){75}}
\put(79,100){\makebox(0,0){$B$}}
\end{picture}

\vspace{0.5cm}

Denote by $m(ABX_1\setminus\Delta)$ the
value of the Lebesgue measure
(on the plane) of $ABX_1\setminus\Delta$.
There are two conceivable cases.
If $m(ABX_1\setminus\Delta) > 0$,
then there is a segment $AA'$ on $l$
from $A$ to a point $A' = (a', 0) \in ABX_1\setminus\Delta$ with
$a < a' < x_n$ for each $x_n$.
There is not any point $X_n$ on $AA'$.

The complement (relative to $l$) of the
set $\{X_n: n \geq 1\}\cup\{A\}$ is a union of
a collection ${\cal I}$ of disjoint open intervals on
$l$. 
Though the lengths of the finite open intervals in ${\cal I}$
converge to 0, no $I \in {\cal I}$ is of length zero.
If $m(ABX_1\setminus\Delta) = 0$,
then there is an open interval $I_a \in {\cal I}$ with
its left endpoint coinciding with $A$ 
(otherwise $m(ABX_1\setminus\Delta) > 0$).
Let $A' = (a', 0)$ be a
point in $I_a$. The segment $AA'$ 
does not contain any point $X_n$.
\proofend

\begin{stat}
\label{st-G2}
There is not a line segment $AA'$
on $l$ from point $A = (a, 0)$ to
point $A' = (a', 0)$, 
where $a' > a$, such that there is not
any point $X_n = (x_n, 0)$ on $AA'$.
\end{stat}

\proof
The point $A$ is the accumulation point of 
the set $\{X_n: n \geq 1\}$ on $l$.
The definition of accumulation point
rules out immediately the existence of such $AA'$.
\proofend

The Statements \ref{st-G1} and 
\ref{st-G2} can also be formulated
and proved in an analytical setting.
Consider a set consisting of terms $x_n, n = 1, 2, \cdots$ of
a strictly decreasing sequence of
positive numbers with limit $a$.
The set $\{x_n: n \geq 1\}\cup\{a\}$ is a closed set. The complement
of $\{x_n: n \geq 1\}\cup\{a\}$ 
is a union of a collection ${\cal I}$ 
of disjoint open intervals.

\begin{stat}
\label{st-A1}
There is an open interval $I_a \in {\cal I}$, such that
the left endpoint of $I_a$ coincides with $a$.
\end{stat}

\proof
Denote by ${\cal J}$ the sub-collection of
all finite open intervals in ${\cal I}$.
Since
\[
\inf\{|a - x|: x \in \cup_{I \in {\cal J}}\} 
= \inf\{|a - x|: x \in I, \; I \in {\cal J}\} 
= 0
\]
there is an $I_a \in {\cal J}$ with 

\[
\inf\{|a - x|: x \in I_a \} = 0.
\]
The lower
endpoint of $I_a$ is $a$.
Actually, such $I_a \in {\cal J}$ with
$\inf I_a = a$ is a consequence of the
result below.

\[
\inf\bigcup_{I \in {\cal J}}I
= \inf\{x: x \in I,\; I \in {\cal J}\} = a.
\]
\proofend

The following statement is an immediate consequence
of the definition of accumulation point.

\begin{stat}
\label{st-A2}
There is not an open interval $I_a \in {\cal I}$, such that
the left endpoint of $I_a$ coincides with $a$.
\end{stat}

A variant of Statement \ref{st-A1} is as follows.

\begin{stat}
\label{st-P1}
There is an open interval $I_a \in {\cal J}$ with a minimum
length.
\end{stat}

We first emphasize the difference between the limit of
a convergent sequence (with different terms)
and any term of the sequence.
Though the lengths of the open intervals in ${\cal J}$
converge to 0, no $I \in {\cal J}$ is of length zero.
Replacing the positive length of any $I \in {\cal J}$
with the limit 0
is strictly prohibited in current mathematics.

\proof
Consider a thought experiment. 
For each open interval
$I \in {\cal J}$, there is a timer with a pre-set time period
equal to the length of $I$.
All timers start at the same instant, and 
also stop simultaneously
once the time period of
any of the timers expires.
Let $\tau$ be the time elapsed between the
starting and stopping instants.
Since $\tau = 0$ implies an open interval 
of a zero length, which contradicts the
Archimedean axiom, we have $\tau > 0$.
There is then 
an open interval $I_a \in {\cal J}$ with a minimum
length $\tau$.
\proofend

\begin{stat}
\label{st-P2}
There is not an open interval $I_a \in {\cal J}$ with a minimum
length.
\end{stat}

\proof
Let $|I|$ represent the length of an interval $I$.
If there is an open interval $I_a \in {\cal J}$ with a minimum
length $\tau$, then

\[
\sum_{I \in {\cal J}}|I| \geq \sum_{I \in {\cal J}}\tau.
\]
The right side is infinite. But
\begin{equation}
\label{eq-sum}
\sum_{I \in {\cal J}}|I| = x_1 - a
\end{equation}
is finite.
\proofend

Current mathematics provides yet
another way to prove Statement \ref{st-P2} 
by an assertion: A countable intersection of 
infinitely many nested open intervals with a common
endpoint is an empty set.
This assertion is also an immediate consequence of
the definition of accumulation point or limit.

Similar to the proofs of Statements \ref{st-P1} and \ref{st-P2}, 
the following two statements are both
deducible. Denote by ${\cal K}$ the collection of
the open intervals such that
$\cup_{I \in {\cal K}}I$
is the complement of the Cantor set relative to
the closed unit interval $[0, 1]$.

\begin{stat}
There is an open interval $I \in {\cal K}$ with a minimum
length.
\end{stat}

\begin{stat}
There is not an open interval $I \in {\cal K}$ with a minimum
length.
\end{stat}

Some statements above
follow immediately from definitions or axioms
adopted in current mathematics. The
derivations of the other statements
rely on some intermediate results, which can
be found in standard textbooks. 
All the statements are logical consequences
of current mathematics, and hence are
mathematically true. 
It is meaningless to use one of the statements 
(or the definitions and axioms behind the statement) to
argue against the other.
The contradictions between the statements should
only be interpreted as
the indication of the inconsistency of
current mathematics.

There might be some psychological
difficulties in
conceiving the inconsistency.
The inconsistency involves
the notion of accumulation point or limit,
which in turn relies on the notion of infinite.
There are two different conceptions
of infinite in mathematics:
potential infinite versus actual infinite.
A difficulty might
originate from potential infinite.

For a convergent sequence of different numbers, the notion of
limit is formulated by the ($N, \epsilon$) argument
(``for each positive real number $\epsilon$, there is
a positive integer $N$ ...'').
Different conceptions of infinite
lead to different interpretations of the
formulation. 
The interpretation of limit based on
potential infinite goes like this.
For each specified positive number $\epsilon$,
a positive integer $N$ can be found with
the required property, and this process of
``finding an $N$ for each specified $\epsilon$''
will never end.
With such interpretation of limit,
one might deny the existence of the collection 
${\cal I}$ (or ${\cal K}$)
of the open intervals in the first
place. But in the eyes of anyone who is trained in
the set-theoretic tradition, such denial might be
absurd, since the existence of ${\cal I}$ (or ${\cal K}$)
seems to be an obvious consequence of the existence of
the set of all real numbers.

On the other hand,
from the viewpoint of actual infinite,
the ($N, \epsilon$) formulation of limit contains,
as remarked by Abraham Robinson,
``a clear reference to a well-defined
infinite totality, i.e., the totality of
positive real numbers'' \cite{Hint}.
However, even for those who agree with
the viewpoint of actual infinite,
there might be yet a difficulty in
conceiving the inconsistency.

For example, 
consider (\ref{eq-sum}).
The right side of (\ref{eq-sum})
is the limit of a sequence of partial sums.
The interpretation of limit based on the
notion of actual infinite, in this example,
treats ${\cal J}$ as a completed totality, and
the limit as the total length of the open intervals
in ${\cal J}$.
So the length of each open interval in ${\cal J}$ 
is a positive addendum
of an infinite sum. The limit is a result of
adding together all the addenda without
leaving any out.
However, is there a last addendum added into
the total? In other words, is there an open interval
in ${\cal J}$ with a minimum length?

By Statement \ref{st-P1}, there is an open interval
in ${\cal J}$ with a minimum length.
As shown in the proof of Statement
$\ref{st-P1}$, the non-existence of
such open interval implies an open interval of length zero,
which contradicts the Archimedean axiom.
But the existence of an open interval in ${\cal J}$
with a minimum length
implies a maximum natural number, which contradicts 
Peano's axioms. 
The basis of
natural numbers and
any reasoning involving natural numbers is a set of assumptions,
called
Peano's axioms. The Archimedean axiom involves natural numbers,
and hence is also a consequence of Peano's axioms.
The inconsistency of current mathematics revealed 
in this paper suggests a fundamental flaw in the
formulation of the notion of natural numbers.

Many mathematicians share a
belief:
All mathematical statements should be reducible ultimately
to statements about natural numbers \cite{Cour}.
Perhaps it is this belief that makes
the inconsistency of current mathematics 
difficult to detect.



\begin{thebibliography}{99} 
\bibitem{Cour}
R. Courant and H. Robbins, revised by
I. Stewart,
{\bf What Is Mathematics?}, 2nd. ed.,
Oxford University Press, Inc., 1996.
\bibitem{Ebbi}
H.-D. Ebbinghaus, J. Flum, and W. Thomas,
{\bf Mathematical Logic}, 2nd ed., Springer-Verlag, 1994.
\bibitem{Hint}
J. Hintikka, ed., 
{\bf The Philosophy of Mathematics}, Oxford, 1969.




\end{thebibliography}
\end{document}